\documentclass[11pt]{article}
\usepackage{enumerate}
\usepackage{amssymb,a4wide,latexsym,makeidx,epsfig,fleqn}
\usepackage{amsthm}
\usepackage{amsmath}
\allowdisplaybreaks[4]
\usepackage{enumerate}
\usepackage{graphicx}
\usepackage{color}
\usepackage{epstopdf}
\newtheorem{theorem}{Theorem}[section]
\newtheorem{remark}[theorem]{Remark}

\newtheorem{lemma}[theorem]{Lemma}

\newtheorem{corollary}[theorem]{Corollary}

\begin{document}
\textwidth 150mm \textheight 230mm
\setlength{\topmargin}{-15mm}
\title{Sufficient conditions for $t$-tough graphs to be Hamiltonian and pancyclic or bipartite
\footnote{This work is supported by the National Natural Science Foundations of China (No. 12371348, 12201258), the Postgraduate Research \& Practice Innovation Program of Jiangsu Province (No. KYCX25\_3155), the Postgraduate Research \& Practice Innovation Program of Jiangsu Normal University (No. 2025XKT0632, 2025XKT0633) .}}
\author{{ Xiangge Liu, Caili Jia, Yong Lu\footnote{Corresponding author}, Jiaxu Zhong }\\
{\small  School of Mathematics and Statistics, Jiangsu Normal University,}\\ {\small  Xuzhou, Jiangsu 221116,
People's Republic
of China.}\\
{\small E-mail:  luyong@jsnu.edu.cn}}

\date{}
\maketitle
\begin{center}
\begin{minipage}{120mm}
\vskip 0.3cm
\begin{center}
{\small {\bf Abstract}}
\end{center}
{\small

The toughness of graph $G$, denoted by $\tau(G)$, is $\tau(G)=\min\{\frac{|S|}{c(G-S)}:S\subseteq V(G),c(G-S)\geq2\}$ for every vertex cut $S$ of $V(G)$ and the number of components of $G$ is denoted by $c(G)$.
Bondy in 1973, suggested the ``metaconjecture"  that almost any nontrivial condition on a graph which implies that the graph is Hamiltonian also implies that the graph is pancyclic.
Recently,  Benediktovich  [Discrete Applied Mathematics. 365 (2025) 130--137] confirmed  the Bondy's metaconjecture for $t$-tough graphs in the case when $t\in\{1;2;3\}$ in terms of the size, the spectral radius and the   signless Laplacian spectral radius of the graph.
In this paper, we will confirm the Bondy's metaconjecture for $t$-tough graphs in the case when  $t\geq4$ in terms of the size, the spectral radius, the signless Laplacian spectral radius, the distance spectral radius and the  distance signless Laplacian  spectral radius of graphs.

 }
\end{minipage}
\end{center}

\section{Introduction }
\emph{}

Let $G=(V(G),E(G))$ be a simple graph, where $V(G)$ is the vertex set and $E(G)$ is the edge set.
The order and size of $G$ are denoted by $|V(G)|=n$ and $|E(G)|=m$, respectively.
For a vertex subset $S$ of $G$, we denote by $G-S$ and $G[S]$ the subgraph of $G$ obtained from $G$ by deleting the vertices in $S$ together with their incident edges and the subgraph of $G$ induced by $S$, respectively.
Let $t$ be a positive real number and a connected graph $G$ is \emph{$t$-tough} if $tc(G-S)\leq|S|$ for every vertex cut $S$ of $V(G)$ and the number of components of $G$ is denoted by $c(G)$.
The \emph{toughness} of graph $G$, denoted by $\tau(G)$, is $\tau(G)=\min\{\frac{|S|}{c(G-S)}:S\subseteq V(G),c(G-S)\geq2\}$ for every vertex cut $S$ of $V(G)$ (taking $\tau(K_{n})=\infty$).
This concept  was initially proposed by Chv\'{a}tal \cite{CV} in 1973, which serves as a simple way to measure how tightly various pieces of a graph hold together.

Let $d(v)$ be the \emph{degree} of vertex $v\in V(G)$, $\delta(G)$ (or $\Delta(G)$) be the \emph{minimum degree} (or \emph{maximum degree}) of $G$ and short for $\delta$ (or $\Delta$).
We denote the \emph{cycle} and the \emph{complete graph} on $n$ vertices by using $C_{n}$ and $K_{n}$, respectively.
For two vertex-disjoint graphs $G_{1}$ and $G_{2}$, we use $G_{1}+G_{2}$ to denote the \emph{disjoint union} of $G_{1}$ and $G_{2}$.
The \emph{join} $G_{1}\vee G_{2}$ is the graph obtained from $G_{1}+G_{2}$ by adding all possible edges between $V(G_{1})$ and $V(G_{2})$.
If $G_{1}=G_{2}=\cdots=G_{k}$, we denote $G_{1}+G_{2}+\cdots+G_{k}$ by $kG_{1}$.
A graph $G$ is called \emph{Hamiltonian} if it has a cycle that contains all vertices of $G$, and \emph{pancyclic} if it has simple cycles of all lengths from 3 to $n$.
A pancyclic graph is certainly Hamiltonian, but not conversely in general.
We denote a \emph{bipartite graph} with bipartition $(X,Y)$ by using $G[X,Y]$, and a \emph{complete bipartite graph} with two parts having $m,n$ vertices by using $K_{m,n}$.

For a simple graph $G$ of order $n$,
the \emph{adjacency matrix} of $G$ is denoted by $A(G)$, and $\lambda_{1}(G)\geq\lambda_{2}(G)\geq\cdots\geq\lambda_{n}(G)$ are the eigenvalues of $A(G)$.
In particular, the eigenvalue $\lambda_{1}(G)$ is called the \emph{spectral radius} of $G$.
The \emph{distance} between $v_{i}$ and $v_{j}$ denoted by $d_{ij}(G)$, is the length of a shortest path from $v_{i}$ to $v_{j}$, and the \emph{distance matrix} of $G$ is the symmetric matrix $D(G)=(d_{ij}(G))_{n\times n}$  with its rows and columns indexed by $V(G)$.
We can order the eigenvalues of $D(G)$ as $\lambda_{1}(D(G))\geq \lambda_{2}(D(G))\geq\cdots\lambda_{n}(D(G))$ where the \emph{distance spectral radius} of $G$ is $\lambda_{1}(D(G))$.
The matrix $Q(G)=D'(G)+A(G)$, where $D'(G)$ is a \emph{degree diagonal matrix} of $G$, is called the \emph{signless Laplacian matrix} of $G$ and the largest eigenvalue of $Q(G)$, denoted by $q(G)$, is called to be the \emph{signless Laplacian spectral radius}.
The \emph{transmission} $Tr(v)$ of a vertex $v\in V(G)$ is the sum of distances from $v$ to all vertices in $G$. We say that a graph is \emph{$k$-transmission-regular} (or transmission-regular) if its distance matrix has constant row sum equal to $k$. The \emph{distance signless Laplacian matrix} of a graph $G$, is given by $Q_{D}(G)=Diag(Tr)+D(G)$,
where $Diag(Tr)$ is the diagonal matrix whose diagonal entries are the vertex transmissions in $G$.   Let $\eta_{1}(G)\geq\eta_{2}(G)\geq\cdots\geq\eta_{n}(G)$ be eigenvalues of $Q_{D}(G)$  and the largest eigenvalue $\eta_{1}(G)$ is also called the \emph{distance signless Laplacian spectral radius}.

Many researchers established the relationship between toughness conditions and the spectra of graphs.
Fan et al. \cite{FLL} presented sufficient conditions based on the spectral radius for  graphs to be  1-tough with minimum degree $\delta(G)$ and graphs to be $t$-tough with $t\geq 1$ being an integer, respectively.
Lou et al. \cite{LLS}  presented a sufficient condition based on the distance spectral radius to guarantee that a graph is $t$-tough, where $t$ or $\frac{1}{t}$ is a positive integer.
Jia and Lou \cite{JL} also considered the above problem  based on the signless Laplacian spectral radius.
For more extensive work on  toughness conditions related to the spectra of graphs, one can see (\cite{BAE},\cite{CLX},\cite{CFL},\cite{CYY}).

Meanwhile, the study of Hamiltonian, pancyclic and the spectra of graphs has received extensive attention.
In 1972, Erd\H{o}s \cite{E} conjectured that every Hamiltonian graph with independence number at most $k$ and at least $n=\Omega(k^{2})$ vertices is pancyclic.
Dragani\'{c} et al. \cite{DCS} proved this  conjecture in a strong form.
Zhou and Wang \cite{ZQW} established some sufficient conditions for a graph to be Hamilton-connected in terms of the size, the spectral radius and the signless Laplacian spectral radius of the graph.
They \cite{ZW}  also considered the relationship between distance signless Laplacian spectral radius and the Hamiltonian properties of graphs according to the size $m$ and the order $n$.
For more extensive work on  Hamiltonian, pancyclic and the spectra of graphs, one can see (\cite{GN},\cite{LSX},\cite{XYA},\cite{ZB},\cite{ZZL}).

Chv\'{a}tal and other scholars investigated the relationship between toughness conditions and the existence of cyclic structures, in particular, determining whether the graph is Hamiltonian or pancyclic.
It is well-known that recognition problems whether the graph is Hamiltonian or pancyclic are NP-complete.
Cai et al. \cite{CYX} established sufficient conditions for a tough graph $G$ to be Hamiltonian  based on the size, the spectral radius and  the  signless Laplacian spectral radius.
Fan et al. \cite{FLL} also presented sufficient conditions for tough graphs to be Hamiltonion in terms of the spectral radius.

As is well known,  if a graph $G$ with  $n\geq3$,  and  every vertex has minimum degree at least $\frac{n}{2}$, then $G$ is Hamiltonian.
Bondy showed this condition implies  that  either $G$ is isomorphic to the complete bipartite graph $K_{\frac{n}{2},\frac{n}{2}}$ or $G$ is pancyclic.
In 1973, Bondy \cite{Bondy} proposed the following ``metaconjecture" which has had a profound effect on research on cycles in graphs:

\noindent\textbf{Bondy's metaconjecture.}
Almost any non-trivial condition on a graph which implies that the graph is Hamiltonian also implies
that the graph is pancyclic.
(There may be a simple family of exceptional graphs.)

Recently,  Benediktovich \cite{BVI} confirmed  the Bondy's metaconjecture  for $t$-tough graphs in the case when $t\in\{1;2;3\}$ in terms of the size, spectral radius and   signless Laplacian spectral radius of the graph.
For more extensive work on  toughness conditions and Hamiltonian or pancyclic, one can see (\cite{BBV},\cite{BPP},\cite{OS},\cite{SH},\cite{SS}).

More generally, motivated the  results of \cite{BVI}, we present sufficient conditions for $t$-tough graphs to be Hamiltonian and pancyclic or bipartite in terms of the size firstly.

\noindent\begin{theorem}\label{th:1.1.}
Let $G$ be a $t$-tough  ($t\geq4$) and simple connected graph with $n>10t-3$ vertices and $m$ edges. If
\begin{align}
m\geq\binom{n-2t}{2}+3t^{2},
\end{align}
then $G$ is Hamiltonian, especially $G$ is  pancyclic or bipartite.
\end{theorem}

Additionally, we establish sufficient conditions based on the spectral radius,  the signless Laplacian spectral radius, the distance spectral radius and the distance signless Laplacian  spectral radius to ensure that $t$-tough graphs to be Hamiltonian and pancyclic or bipartite, respectively.

\noindent\begin{theorem}\label{th:1.2.}
Let $G$ be a $t$-tough graph ($t\geq4$) and simple connected with $n>10t-3$ vertices and $m$ edges. If
\begin{enumerate}
\item [(i)] $\lambda_{1}(G)\geq\sqrt{n^{2}-(4t+2)n+10t^{2}+2t+1},$ or
\item [(ii)] $q(G)\geq\frac{n^{2}-(4t+1)n+10t^{2}+2t}{n-1}+n-2,$ or
\item [(iii)] $\lambda_{1}(D(G))\leq\ n+4t-1-\frac{10t^{2}+2t}{n},$ or
\item [(iv)] $\eta_{1}(G)\leq\ 2n+8t-2-\frac{20t^{2}+4t}{n},$
\end{enumerate}
then $G$ is Hamiltonian, especially $G$ is pancyclic or bipartite.
\end{theorem}

\section{Proof of Theorem \ref{th:1.1.} }
In this section, we prove sufficient condition for $t$-tough graph to be Hamiltonian and pancyclic or bipartite  in terms of the size. At first, we will give some definitions and lemmas.

For a nonnegative integer $k$, the \emph{$k$-closure} of a graph $G$, denote by $\mathcal{C}_{k}(G)$, is obtained from $G$ by recursively joining pairs of nonadjacent vertices the degree sum of which is at least $k$ as far as such a pair remains.
Bondy and Chv\'{a}tal showed that $\mathcal{C}_{k}(G)$ is uniquely defined by a given graph $G$.

\noindent\begin{lemma}\label{le:2.1.}\cite{BV}
A graph $G$ is Hamiltonian if and only if $\mathcal{C}_{n}(G)$ is Hamiltonian.
\end{lemma}

Let $(d_{1},d_{2},\ldots,d_{n})$ be a nondecreasing degree sequence of $G$ that is, $d_{1}\leq d_{2}\leq\cdots\leq d_{n}$. For convenience, we use $(0^{x_{0}},1^{x_{1}},\ldots,k^{x_{k}},\ldots,\Delta^{x_{\Delta}})$ to denote the degree sequence of $G$, where $x_{k}$ is the number of vertices of degree $k$ in the graph $G$.

The \emph{predicate} $P(t)$ was defined as
$P(t)$: for all $k<\frac{n}{2}$, if $d_{k}\leq k$ then $d_{n-k+t}\geq n-k$.

\noindent\begin{lemma}\label{le:2.2.}\cite{SA}
Let $t\geq4$ be a positive integer.  If $G$ is a $t$-tough graph satisfying $P(t)$, then $G$ is Hamiltonian.
\end{lemma}

\noindent\begin{lemma}\label{le:2.3.}\cite{H}
Let $t$ be a positive integer. If a $t$-tough graph $G$ satisfies $P(t)$ and is Hamiltonian, then $G$ is pancyclic or bipartite.
\end{lemma}

\noindent\begin{lemma}\label{le:2.4.}\cite{H}
If $G$ is not complete $t$-tough graph with a minimum degree $\delta(G)$, then $\delta(G)\geq2t$.
\end{lemma}

\noindent\begin{lemma}\label{le:2.5.}\cite{J}
Let $G$ be a graph without a Hamiltonian cycle and at least 11 vertices. Then
\begin{enumerate}
\item there exist two non-adjacent vertices $x$, $y$ such that $d(x)+d(y)\leq\mid V(G)\mid-5$, or
\item there exist for some $t\geq1$ vertices $x_{1},x_{2},\ldots,x_{t}$ such that $G-x_{1}-\cdots-x_{t}$ has at least $t+1$ components.
\end{enumerate}
\end{lemma}

According to Lemma \ref{le:2.5.}, we obtain the following corollary, which will play an essential role in the proof of Theorem \ref{th:1.1.}.

\noindent\begin{corollary}\label{co:2.6.}
Let $G$ be a $t$-tough graph without a Hamiltonian cycle with at least 11 vertices. Then there exist two non-adjacent vertices $x$, $y$ such that $d(x)+d(y)\leq\mid V(G)\mid-5$.
\end{corollary}

\noindent\textbf{Proof.}
Since $G$ has no Hamiltonian cycle, $G$ is not a complete graph $K_{n}$. If there exist for some $s\geq1$ vertices $x_{1},x_{2},\ldots,x_{s}$ such that $G-x_{1}-\cdots-x_{s}$ has at least $s+1$ components,
by
$$\tau(G)=\min\{\frac{\mid S\mid}{c(G-S)}:S\subseteq V(G),c(G-S)\geq2\},$$
then
$$t\leq\tau(G)\leq\frac{s}{s+1}<1,$$
a contradiction.

The proof is completed.
$\hfill\square$\\

\noindent\begin{lemma}\label{le:2.7.} \cite{H}
Let Hamiltonian graph $G$ of order $n$ contain more than $\frac{n}{3}$ vertices of degree greater than $\frac{n}{2}$. Then graph $G$ is pancyclic.
\end{lemma}

\noindent\begin{lemma}\label{le:2.8.} \cite{BVI}
2-tough graph $G$ is Hamiltonian if and only if its $(n-1)$-closure $\mathcal{C}_{n-1}(G)$ is Hamiltonian.
\end{lemma}

Now, we will prove Theorem \ref{th:1.1.}.\\

\noindent\textbf{Proof of Theorem \ref{th:1.1.}.}

Suppose to the contrary, $G$ is neither  Hamiltonian nor pancyclic and bipartite.
Thus it does not satisfy the predicate $P(t)$, since otherwise by Lemmas \ref{le:2.2.} and  \ref{le:2.3.}, it would be Hamiltonian, pancyclic or bipartite.
Hence, there is a positive integer $k$, such that $k<\frac{n}{2}$, $d_{k}\leq k$ and $d_{n-k+t}\leq n-k-1$.
Therefore, we have
\begin{align}
2m&=\sum_{i=1}^{k}d_{i}+\sum_{i=k+1}^{n-k+t}d_{i}+\sum_{i=n-k+t+1}^{n}d_{i}\nonumber
\\&\leq k^{2}+(n-2k+t)(n-k-1)+(k-t)(n-1)\nonumber
\\&=n^{2}-n+3k^{2}+(1-2n-t)k\nonumber
\\&=2\binom{n-2t}{2}+6t^{2}-(k-2t)(2n-3k-5t-1).
\nonumber
\end{align}

Thus
\begin{align} m\leq\binom{n-2t}{2}+3t^{2}-\frac{(k-2t)(2n-3k-5t-1)}{2}.
\end{align}

Combining with the condition of Theorem \ref{th:1.1.}, we know that $\binom{n-2t}{2}+3t^{2}\leq m\leq\binom{n-2t}{2}+3t^{2}-\frac{(k-2t)(2n-3k-5t-1)}{2}$, thus $(k-2t)(2n-3k-5t-1)\leq0$.
By Lemma \ref{le:2.4.}, we have $\delta(G)\geq2t$.
Since $\delta(G)\leq d_{k}\leq k$, we obtain $k\geq2t$.\\

\textbf{Case 1.} $(k-2t)(2n-3k-5t-1)=0$.

In this case, we have $k=2t$, or $k\ne2t$ and $(2n-3k-5t-1)=0$.

\textbf{Subcase 1.1.} $k=2t$.

If $k=2t$, then $G$ is a graph with $d_{2t}\leq 2t$, $d_{n-t}\leq n-2t-1$ and $d_{n}\leq n-1$. Hence, inequalities in (1) and (2)
which have the form $\binom{n-2t}{2}+3t^{2}\leq m\leq\binom{n-2t}{2}+3t^{2}$, are equalities.
Therefore, $G$ has degree sequence $((2t)^{2t},(n-2t-1)^{n-3t},(n-1)^{t})$.

Let us show that $(n-1)$-closure $\mathcal{C}_{n-1}(G)$ of graph $G$ is the complete graph $\mathcal{C}_{n-1}(G)={K}_{n}$. Since $(n-2t-1)+(n-2t-1)=2n-4t-2>n-1$, all $(n-3t)$ vertices of degree $(n-2t-1)$ must be pairwise adjacent in $\mathcal{C}_{n-1}(G)$.
Then each of these vertices must be adjacent to each of $2t$ vertices of degree $2t$ in $\mathcal{C}_{n-1}(G)$. Hence, these $2t$ vertices have degree at least $(n-3t)+t=n-2t$ in $\mathcal{C}_{n-1}(G)$.
Therefore, each of these $2t$ vertices must be pairwise adjacent in $\mathcal{C}_{n-1}(G)$, that is $\mathcal{C}_{n-1}(G)={K}_{n}$ which implies Hamiltonian.
Since  $G$ is $t$-tough ($t\geq4$), then it is 2-tough.
By Lemmas \ref{le:2.7.}  and \ref{le:2.8.}, $G$ is Hamiltonian, and then  $G$ is pancyclic, a contradiction.

\textbf{Subcase 1.2.} $k\ne2t$ and~$(2n-3k-5t-1)=0$.

Since $k<\frac{n}{2}$ and $n>10t-3$, we can get $4k+2\leq 2n=3k+5t+1$ and $20t-6<2n=3k+5t+1$, hence $5t-2\leq k\leq 5t-1$.
Due to $n=\frac{3k+5t+1}{2}$ is a positive integer, we obtain $k=5t-1$ and $n=10t-1$.
Then $d_{5t-1}\leq 5t-1$, $d_{6t}\leq 5t-1$ and  $d_{10t-1}\leq 10t-2$, and all inequalities of $(2)$ become equalities.
Thus  $G$ has degree sequence $((5t-1)^{6t},(10t-2)^{4t-1})$.
Since there are no two vertices $x,y$ that are non-adjacent such that $d(x)+d(y)\leq |V(G)|-5$.
By Corollary \ref{co:2.6.} and Lemma \ref{le:2.7.} , we can get $G$ is Hamiltonian, and then  $G$ is pancyclic, a contradiction.\\

\textbf{Case 2.} $(k-2t)(2n-3k-5t-1)<0$.

In this case, $k\geq 2t$ and~$(2n-3k-5t-1)<0$.
Since $k<\frac{n}{2}$, we have $n\geq2k+1$.
According to  these results, we have $4k+2\leq2n\leq3k+5t$, thus $k\leq5t-2$ and $n\leq10t-3$, a contradiction.

Combining with Case 1 and Case 2, a $t$-tough graph with $n>10t-3$ vertices, if $m\geq\binom{n-2t}{2}+3t^{2}$, then $G$ is Hamiltonian, especially $G$ is pancyclic or bipartite.

This completes the proof.
$\hfill\square$\\

In particular, if $t\in\{1;2;3\}$, then we can obtain the theorem as follows which confirmed by  Benediktovich in \cite{BVI}.

\textbf{[Theorem 2.1 in \cite{BVI}]}
Let $G$ be a $t$-tough  and simple connected graph with $n$ vertices and $m$ edges. If $$m\geq\binom{n-2t}{2}+3t^{2},$$ then  $G$ is pancyclic or bipartite
\begin{enumerate}
    \item[(1)]{when $t=1$ and $n\geq 7$;}
    \item[(2)]{when $t=2$ and $n\geq 16$;}
    \item[(3)]{when $t=3$ and $n\geq28$.}
\end{enumerate}

In this theorem,   Benediktovich discussed the special graphical degree sequences when $t=1$ and $n=7$, $t=2$ and $n=16,17$ which can be verified  pancyclic  or bipartite. Note that, in Theorem \ref{th:1.1.}, although we restrict $t\geq4$, when $t=3$ and $n>27$, i.e., $n\geq28$, our results coincide with those in \cite{BVI}.   \\

\noindent\begin{remark}
In order to verify our theorem and show our conclusions better, we use another method to prove the sufficient conditions for 4-tough graphs to be Hamiltonian, pancyclic or bipartite, and we display the corresponding graphs (see Figure 1) during the proof.\\
\end{remark}

\begin{figure}[htbp]
  \centering
  \includegraphics[scale=1]{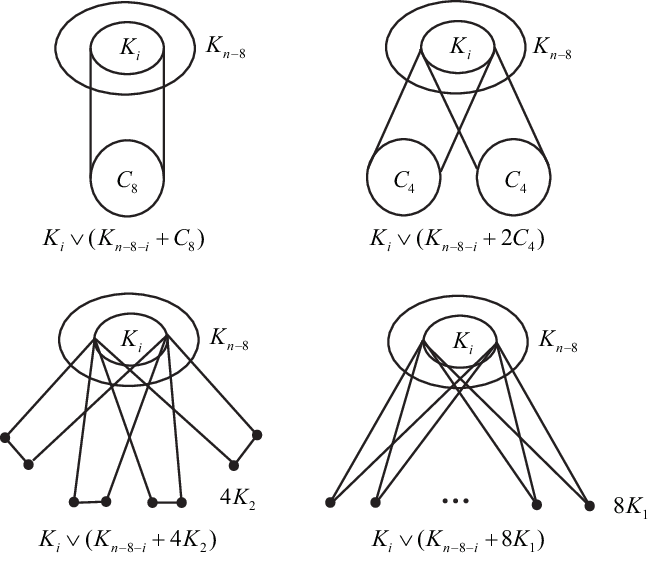}
  \caption{Corresponding graphs in the proof of $t=4$.}
\end{figure}

\noindent\textbf{Proof.}

Suppose to the contrary, $G$ is not Hamiltonian and has degree sequence $d_{1}\leq d_{2}\leq \cdots\leq d_{n}$.
According to the proof of Theorem \ref{th:1.1.}, we get $(k-2t)(2n-3k-5t-1)\leq0$.
By Lemma \ref{le:2.4.}, $\delta(G)\geq2t=8$.
Since $\delta(G)\leq d_{k}\leq k$, we obtain $k\geq8$.
Thus, we divide the following proof into two cases.\\

\textbf{Case 1.}
$(k-2t)(2n-3k-5t-1)=0$.

In this case, we have $k=8$ or $k\neq8$ and $2n-3k-21=0$.

\textbf{Subcase 1.1.}
$k=8$.

If $k=8$, then $G$ is a graph with $d_{8}\leq8$, $d_{n-4}\leq n-9$, $d_{n}\leq n-1$. We have $\binom{n-2t}{2}+3t^{2}\leq m\leq\binom{n-2t}{2}+3t^{2}$ by (1) and (2). Thus all inequalities become equalities and  $G$ has degree sequence $(8^{8},(n-9)^{n-12},(n-1)^{4})$.

Let $S$ be the set containing eight 8-degree vertices of $G$, and there exist two non-adjacent vertices in $S$ by Corollary \ref{co:2.6.}.

By Lemma \ref{le:2.1.} and $G$ is not Hamiltonian, we have $\mathcal{C}_{n}(G)$ is also not Hamiltonian.
Due to the definition of $\mathcal{C}_{n}(G)$, all points except the vertices of $S$ form a complete graph $K_{n-8}$.
If one vertex $v_{i}$ of the $K_{n-8}$ is adjacent to one vertex of $S$, then it must be adjacent to all vertices of $S$.
Moreover, there are at least four vertices of the $K_{n-8}$ adjacent to all of $S$.

If there are four vertices of the $K_{n-8}$ adjacent to all vertices of $S$, then $\mathcal{C}_{n}(G)\supseteq K_{4}\vee(K_{n-12}+ C_{8})$ or $K_{4}\vee(K_{n-12}+ 2C_{4})$. Since $K_{4}\vee(K_{n-12}+ C_{8})$ and $K_{4}\vee(K_{n-12}+ 2C_{4})$ are Hamiltonian, then $\mathcal{C}_{n}(G)$ is Hamiltonian, a contradiction.

If there are five vertices of the $K_{n-8}$ adjacent to all vertices of $S$, then $\mathcal{C}_{n}(G)\supseteq K_{5}\vee(K_{n-13}+ C_{8})$ or $K_{5}\vee(K_{n-13}+ 2C_{4})$ or $K_{5}\vee(K_{n-13}+ 4K_{2})$. Since $K_{5}\vee(K_{n-13}+ C_{8})$, $K_{5}\vee(K_{n-13}+ 2C_{4})$ and $K_{5}\vee(K_{n-13}+ 4K_{2})$ are Hamiltonian, then $\mathcal{C}_{n}(G)$ is Hamiltonian, a contradiction.

If there are six vertices of the $K_{n-8}$ adjacent to all vertices of $S$, then $\mathcal{C}_{n}(G)\supseteq K_{6}\vee(K_{n-14}+ C_{8})$ or $K_{6}\vee(K_{n-14}+ 2C_{4})$ or $K_{6}\vee(K_{n-14}+ 4K_{2})$.
Since $K_{6}\vee(K_{n-14}+ C_{8})$, $K_{6}\vee(K_{n-14}+ 2C_{4})$ and $K_{6}\vee(K_{n-14}+ 4K_{2})$ are Hamiltonian, then $\mathcal{C}_{n}(G)$ is Hamiltonian, a contradiction.

If there are seven vertices of the $K_{n-8}$ adjacent to all vertices of $S$, then $\mathcal{C}_{n}(G)\supseteq K_{7}\vee(K_{n-15}+ 4K_{2})$. Since $K_{7}\vee(K_{n-15}+ 4K_{2})$ is Hamiltonian, then $\mathcal{C}_{n}(G)$ is Hamiltonian, a contradiction.

If there are eight vertices of the $K_{n-8}$ adjacent to all vertices of $S$, then $\mathcal{C}_{n}(G)\supseteq K_{8}\vee(K_{n-16}+ 8K_{1})$. Since $K_{8}\vee(K_{n-16}+ 8K_{1})$ is Hamiltonian, then $\mathcal{C}_{n}(G)$ is Hamiltonian, a contradiction.

If there are more than eight vertices of $K_{n-8}$ adjacent to all vertices,  then $\mathcal{C}_{n}(G)\supseteq K_{i}\vee(K_{n-8-i}+ 8K_{1})$ ($i>8$). Since $K_{i}\vee(K_{n-8-i}+ 8K_{1})$ ($i>8$) is Hamiltonian, then $\mathcal{C}_{n}(G)$ is Hamiltonian, a contradiction.

\textbf{Subcase 1.2.}
$k\neq0$ and $2n-3k-21=0$.

In this case, since $k=\frac{2n-21}{3}<\frac{n}{2}$ and $n>37$, we can get $38\leq n\leq41$.
Thus, let $n=39$ and $k=19$, then $d_{19}\leq19$, $d_{24}\leq19$ and $d_{39}\leq38$.
We have $1026\leq2m=\sum\limits_{i=1}^{39}d_{i}\leq1026$ by (1) and (2), so $\sum\limits_{i=1}^{39}d_{i}=1026$.
The corresponding permissible degree sequence is $(19^{24},38^{15})$.
For this degree sequence, we know that there are no two vertices $x,y$ that are non-adjacent such that $d(x)+d(y)\leq |V(G)|-5$.
By Corollary \ref{co:2.6.}, we can get $G$ is Hamiltonian, a contraction.\\

\textbf{Case 2.}
$(k-8)(2n-3k-21)<0$.

In this case, $k\geq9$ and $2n-3k-21<0$.
Since $k<\frac{n}{2}$, we have $n\geq2k+1$.
From these results, we have $4k+2\leq2n\leq3k+20$, thus $k\leq18$ and $n\leq37$, a contradiction.

Combining with Case 1 and Case 2, a 4-tough graph with $n>37$ vertices, if $m\geq\binom{n-2t}{2}+3t^{2}$, then $G$ is Hamiltonian.

By Lemma \ref{le:2.3.} and  $G$ is Hamiltonian, we obtain $G$ is  pancyclic or bipartite.

This completes the proof.

$\hfill\square$\\

\section{Proof of Theorem \ref{th:1.2.} }

In order to prove Theorem \ref{th:1.2.}, we list some bounds with regard to the spectral radius, the signless Laplacian spectral radius, the distance spectral radius and the  distance signless Laplacian  spectral radius of graphs, respectively.

\noindent\begin{lemma}\label{le:3.1.}\cite{YH}
Let $G$ be a connected graph with $n$ vertices and $m$ edges. Then $\lambda_{1}(G)\leq\sqrt{2m-n+1}$, and the equality holds if and only if $G=K_{n}$ or $G=K_{1,n-1}$.
\end{lemma}

\noindent\begin{lemma}\label{le:3.2.}\cite{YY}
Let $G$ be a graph with $n$ vertices and $m$ edges, $v\in V(G)$. Then $q(G)\leq\frac{2m}{n-1}+n-2$. If $G$ is connected, the equality holds if and only if $G=K_{n}$ or $G=K_{1,n-1}$. Otherwise, the equality holds if and only if $G=K_{n-1}+v$.
\end{lemma}

The \emph{Wiener index} $W(G)$ of a connected graph $G$ of order $n$ is defined by the sum of all distances in $G$, that is, $W(G)=\sum\limits_{i<j}d_{ij}(G)$.

\noindent\begin{lemma}\label{le:3.3.}\cite{IG}
Let $G$ be a connected graph with order $n$. Then $\lambda_{1}(D(G)) \geq \frac{2W(G)}{n}$.
\end{lemma}

\noindent\begin{lemma}\label{le:3.4.}\cite{XZL}
Let $G$ be a connected graph with order $n$. Then
$\eta_{1}(G) \geq \frac{4W(G)}{n}$
with equality if and only if G is transmission regular.
\end{lemma}

Now, we will prove Theorem \ref{th:1.2.}.\\

\noindent\textbf{Proof of Theorem \ref{th:1.2.}.(i)}

By  Lemma \ref{le:3.1.}, the graph $G$ is neither graph $K_{n}$ nor $K_{1,n-1}$.
Combining with the condition of Theorem \ref{th:1.2.}, we have inequalities
$$\sqrt{n^{2}-(4t+2)n+10t^{2}+2t+1}\leq \lambda_{1}(G)<\sqrt{2m-n+1},$$
which implies  $m>\binom{n-2t}{2}+3t^{2}$.
 According to Theorem \ref{th:1.1.}, $G$ is Hamiltonian, and then $G$ is pancyclic or bipartite.
$\hfill\square$\\

In particular, if $t\in\{1;2;3\}$, then we can obtain the theorem as follows which confirmed by  Benediktovich in \cite{BVI}.

\textbf{[Theorem 2.2 in \cite{BVI}]}
Let $G$ be a $t$-tough  and simple connected graph with $n$ vertices and $m$ edges. If $$\lambda_{1}(G)\geq\sqrt{n^{2}-(4t+2)n+10t^{2}+2t+1},$$  then  $G$ is pancyclic or bipartite
\begin{enumerate}
    \item[(1)]{when $t=1$ and $n\geq 7$;}
    \item[(2)]{when $t=2$ and $n\geq 16$;}
    \item[(3)]{when $t=3$ and $n\geq28$.}
\end{enumerate}

In this theorem,   Benediktovich discussed the special graphical degree sequences when $t=1$ and $n=7$, $t=2$ and $n=16,17$ which can be verified pancyclic or bipartite.\\

\noindent\textbf{Proof of Theorem \ref{th:1.2.}.(ii)}

By Lemma \ref{le:3.2.},  $G$ is neither graph $K_{n}$ nor $K_{1,n-1}$.
Combining with the condition of Theorem \ref{th:1.2.}, we have inequalities
$$\frac{n^{2}-(4t+1)n+10t^{2}+2t}{n-1}+n-2\leq q(G)<\frac{2m}{n-1}+n-2,$$
which implies  $m>\binom{n-2t}{2}+3t^{2}$.
According to Theorem \ref{th:1.1.},  $G$ is Hamiltonian, and then $G$ is pancyclic or bipartite.
$\hfill\square$\\

In particular, if $t\in\{1;2;3\}$, then we can obtain the theorem as follows which confirmed by  Benediktovich in \cite{BVI}.

\textbf{[Theorem 2.3 in \cite{BVI}]}
Let $G$ be a $t$-tough  and simple connected graph with $n$ vertices and $m$ edges. If $$q(G)\geq\frac{n^{2}-(4t+1)n+10t^{2}+2t}{n-1}+n-2,$$ then   $G$ is pancyclic or bipartite
\begin{enumerate}
    \item[(1)]{when $t=1$ and $n\geq 7$;}
    \item[(2)]{when $t=2$ and $n\geq 16$;}
    \item[(3)]{when $t=3$ and $n\geq28$.}
\end{enumerate}

In this theorem,   Benediktovich discussed the special graphical degree sequences when $t=1$ and $n=7$, $t=2$ and $n=16,17$ which can be verified pancyclic or bipartite.  It should be noted that in the condition of Theorem 2.3 in \cite{BVI}, the author wrote $$q(G)\geq\frac{2n^{2}-(4t+1)n+10t^{2}+2t}{n-1}+n-2,$$ but it should be $$q(G)\geq\frac{n^{2}-(4t+1)n+10t^{2}+2t}{n-1}+n-2.$$

Note that, in Theorem \ref{th:1.2.}, although we restrict $t\geq4$, when $t=3$ and $n>27$, i.e., $n\geq28$, our results coincide with those in \cite{BVI}.\\

\noindent\textbf{Proof of Theorem \ref{th:1.2.}.(iii)}

Since $Tr(v)\geq d(v)+2(n-1-d(v))=2(n-1)-d(v)$ with equality holds if and only if the maximum distance between $v$ and other vertices in $G$ is at most 2.

So
$$W(G)=\frac{1}{2}\sum\limits_{v\in V(G)}Tr(v)\geq\frac{1}{2}\sum\limits_{v\in V(G)}[2(n-1)-d(v)]=n(n-1)-m$$
with equality holds if and only if the maximum distance between $v$ and other vertices in $G$ is at most 2.

By Lemma \ref{le:3.3.}, we get $\lambda_{1}(D(G))\geq\frac{2W(G)}{n}\geq2(n-1)-\frac{2m}{n}$.
Combining with the condition of Theorem \ref{th:1.2.}, we have inequalities
$$2(n-1)-\frac{2m}{n}\le\lambda_{1}(D(G))\leq n+4t-1-\frac{10t^{2}+2t}{n},$$
which implies $m\geq\binom{n-2t}{2}+3t^{2}$.
According to Theorem \ref{th:1.1.}, $G$ is  Hamiltonian, and then $G$ is pancyclic or bipartite.
$\hfill\square$\\

\noindent\textbf{Proof of Theorem \ref{th:1.2.}.(iv)}

By the above analyses and  Lemma \ref{le:3.4.}, we have $\eta_{1}(G)\geq\frac{4W(G)}{n}\geq4(n-1)-\frac{4m}{n}$.
Combining with the condition of Theorem \ref{th:1.2.},  we have inequalities
$$4(n-1)-\frac{4m}{n}\le\eta_{1}(G)\le 2n+8t-2-\frac{20t^{2}+4t}{n},$$
which implies  $m\ge\binom{n-2t}{2}+3t^{2}$.
According to  Theorem \ref{th:1.1.}, $G$ is  Hamiltonian, and then $G$ is pancyclic or bipartite.
$\hfill\square$\\

\end{document}